\documentclass[reqno, 11pt]{amsart}

\newtheorem{theorem}{Theorem}[section]
\newtheorem{corollary}[theorem]{Corollary}
\newtheorem{lemma}[theorem]{Lemma}

\theoremstyle{definition}
\newtheorem{remark}[theorem]{Remark}

\theoremstyle{definition}
\newtheorem{definition}[theorem]{Definition}

\theoremstyle{definition}
\newtheorem{assumption}[theorem]{Assumption}

\numberwithin{equation}{section}

\newcommand{\Var}{\mathop{\hbox{\rm Var\,}}}
\newcommand\loc{\rm loc}

\newcommand{\nlimsup}{\operatornamewithlimits{\overline{lim}}}

 \makeatletter
 \def\dashint{%
 \operatorname%
 {\,\,\text{\bf--}\kern-.98em\DOTSI\intop\ilimits@\!\!}}
 \makeatother

\newcommand\tr{\hbox{\rm tr}_{2}\,}

\def\sfa{{\sf a}}

\def\bM{\mathbb{M}}
\def\bQ{\mathbb{Q}}
\def\bR{\mathbb{R}}

\def\cF{\mathcal{F}}

\def\cN{\mathcal{N}}

\title[On diffusions   and ``good'' solutions
of elliptic
equations]
{On diffusion processes with
$B(\bR^{2}, VMO)$   coefficients 
and ``good'' Green's functions of the corresponding operators}

\author{N.V. Krylov}

\address{127 Vincent Hall, University of Minnesota, Minneapolis, MN
 55455}
\email{krylov@math.umn.edu}



\begin{document}

\begin{abstract}
The solvability in Sobolev spaces 
with special mixed norms is proved for
nondivergence form second order parabolic equations.
 The leading coefficients are assumed to be measurable in the
time variable and two coordinates of space variables, and 
be almost in
VMO (vanishing mean oscillation) with respect to the other
coordinates. This solvability result
implies the weak uniqueness of solutions
of the corresponding stochastic It\^o equations in the class
of ``good'' solutions   (which is nonempty).
This also implies uniqueness of a Green's
function in the class of ``good'' ones (which is always nonempty).
\end{abstract}

\maketitle

\section{Introduction}
                                                \label{sec1}

Let   $\bR^{d}$ be a Euclidean space of points
 $x=(x^{1},...,x^{d})$ and let $d\geq3$.  
We write
$x=(x',x'')$, where $x'=(x^{1},x^{2})$ and
$x''=(x^{3},...,x^{d})$.
Set
$$
D_{i}u=u_{x^i},\quad D_{ij }u=u_{x^ix^j},\quad
 \partial_{t}u =\partial u/\partial t.
$$
By $Du$ and $D^{2}u$ we
mean the gradient and the Hessian matrix
of $u$.

In this paper    we are dealing with
diffusion processes corresponding to
parabolic  equations in nondivergence form:
\begin{equation}
                                                \label{parabolic}
L u-\lambda u=f,
\end{equation}
where $\lambda\ge 0$ is a constant, 
$f\in L_{p,q} $ (space defined later),  and
\begin{equation}
                                                    \label{9.30.3}
L u= \partial_{t}u +a^{ij }D_{ij }u+b^{i}D_iu-cu.
\end{equation}
We assume that
 all the coefficients are   measurable  and  
$$
c\geq0,\quad |b^i|+ c \le  K,
\quad a^{ij }=a^{ ji},
\quad\delta|\xi|^{2}\leq
a^{rs }\xi^{r}\xi^{s}\leq\delta^{-1}|\xi|^{2} 
$$
for all $i,j=1,...,d$, $\xi\in\bR^{d}$,
where $K$ and $\delta>0$ are fixed constants. 

For  $p,q\in(1,\infty)$
we introduce $L_{p,q} $ as the space of (measurable) functions on $\bR^{d+1}
=\{(t,x):t\in\bR, x\in\bR^{d}\}$
with finite norm given by
$$
\|u\|_{L_{p,q} }^{q}
=\int_{\bR^{d-2}}\Big( \int_{\bR^{3}}|u|^{p}\,dx'dt\big)^{q/p}\,dx''.
$$

Then we introduce the function space in which
we are going to consider $L$ by setting  
$$
W_{p,q}^{1,2} =
\{u:\,u,\partial_{t}u,Du,D^2u\in L_{p,q} \}.
$$

One of the main motivations to consider these particular
Sobolev spaces with
mixed norm  comes from the theory of
stochastic diffusion processes.
Namely, we know that for any $(t,x)\in\bR^{d+1}$
there exists a probability space $(\Omega,\cF,P)$, a
$d$-dimensional random process $x_{s},s\geq0$,
and a
$d$-dimensional Wiener process $w_{s},s\geq0$,
such that $w_{s+h}-w_{s}$ is independent
of $\cF_{s}:=\sigma\{x_{r};r\leq s\}$ for $s,h\geq0$,
$w_{s}$ is $\cF_{s}$-measurable, and, with
probability one, for all $s\geq0$
\begin{equation}
                                                \label{1.7.1}
x_{s}=x+\int_{0}^{s}a^{1/2}(t+r,x_{r})\,dw_{r}
+\int_{0}^{s} b(t+r,x_{r})\,dr,
\end{equation}
where $a=(a^{ij})$, $b=(b^{i})$.

We consider the case that  $a $ is  
only measurable in $(t,x')$. It is 
more regular with respect to $x''$, it 
is almost in VMO.
In such situations, as Uraltseva's examples show, it is not possible to prove
the solvability of \eqref{parabolic} in usual
Sobolev spaces $W^{1,2}_{p}(\bR^{d+1})$ with $p>2+\gamma$,
where $\gamma>0$ is independent of $\delta$. Therefore,
the usual method of proving the weak uniqueness
for \eqref{1.7.1} based on applying It\^o's formula
to $u(t+s,x_{s})$, where $u$ is a solution of \eqref{parabolic},
cannot be justified if $d\geq3$. The common knowledge
until now is that $u$ should be
in $W^{1,2}_{p}$ with $p<d+1$ but sufficiently  
 close to $d+1$ in order to apply It\^o's
formula to $u(t+s,x_{s})$.

 At the same time it is proved in \cite{DK_10} that
\eqref{parabolic} is uniquely solvable in 
$W^{1,2}_{p}(\bR^{d+1})$ with $2<p<2+\gamma$,
where $\gamma=\gamma(d,\delta)>0$,
for any $f\in L_{p}(\bR^{d+1})$. Then a very natural
and puzzling question arose: can equation
\eqref{1.7.1} on perhaps different probability spaces 
have solutions with different distributions
when  the corresponding 
parabolic equation is   uniquely solvable
in some space containing $C^{\infty}_{0}(\bR^{d+1})$. 

Under our assumptions we still do not know the answer
to this question. We only show that weak uniqueness
for \eqref{1.7.1} holds in the set of solutions
whose Green's functions belongs to
$L_{p',q'}$, where $p'=p/(p-1),q'=q/(q-1)$
and $p$ and $q$ will be specified later.
 Of course, we also prove that under our assumptions
 such solutions
do exist and, as a consequence, if for any reason
weak uniqueness holds for \eqref{1.7.1}, then
its Green's function is in $L_{p',q'}$.
This consequence is new even for the equations
with the coefficients measurable in time
and VMO in space variables and, for that matter,
new in the case that
the coefficient are measurable in time
and  continuous  in $x$ uniformly
with respect to $t$, the classical case
treated in \cite{SV_79}.
Analogously, we introduce the notion of ``good''
Green's functions for the operator $L$
and prove that such ``good'' functions exist
and are unique.

Our methods are based on \cite{DK_10} and 
a simple consequence of the
Rubio
de Francia extrapolation theorem
presented in \cite{DK_18} or \cite{DK_18_1}.

\section{Main results}
                                \label{mainsec}

First we introduce some notation.  On many occasions we need to
 take derivatives with respect to  only
part of variables. The reader understands
the meaning of the following notation:
$$
D_{x'}u=u_{x'},\quad D_{x''}u=u_{x''},\quad
D^{2}_{x' }u =D_{x'x'}u=u_{x'x'},
$$
$$
D_{x'x''}u=u_{x'x''},\quad
D^{2}_{ x'' }u=D_{x''x''}u=u_{x''x''}.
$$

$$
W_{p}^{1,2}(\bR^{d+1})=
\{u:\,u,Du,D^2u,\partial_{t}u\in L_{p}(\bR^{d+1})\}.
$$
We also use the abbreviations
$$
C^{\infty}_{0}=C^{\infty}_{0}(\bR^{d+1}),\quad
L_{p}=L_{p}(\bR^{d+1}),\quad
W^{1,2}_p=W^{1,2}_p(\bR^{d+1}),...
$$
For   matrix-valued functions $a(t,x)$ on $\bR^{d+1}$ we
understand $\|a\|_{L_{p}}^{p}$ as
$$
\int_{\bR^{d+1}}|\mbox{trace}\,aa^{*}|^{p/2}\,dx\,dt.
$$
Accordingly are introduced the norms in $W$ spaces.

If $B$ is a Borel subset of a  plane $\Gamma$ in a Euclidean space,
we denote by $|B|$ its volume relative to $\Gamma$. This notation is somewhat
ambiguous because $B$ also belongs to the ambient space,
where its volume can be zero. However, we hope that
from the context it will be clear relative to which  plane
we take the volume in each instance. If there is a measurable
function $f$ on $B$ which is integrable with respect to
the Lebesgue measure $\ell$ on $\Gamma$ we set
$$
 f _{B} =\dashint_{B}f( x)\,\ell(dx):=\frac{1}{|B|}
\int_{B}f( x)\,\ell(dx).
$$

 Let
$$
B_r'(x') = \{ y' \in \bR^{2}: |x' -y' | < r\},
$$
$$
B_r''(x'') = \{ y'' \in \bR^{d-2}: |x''-y''| < r\},\quad B_{r}(x)=
B_r'(x')\times B_r''(x''),
$$
$$
Q_r(t,x) = (t+r^2,t) \times B_r(x),
\quad Q_r=Q_r(0,0),
$$
and let $\bQ$ be the collection of all $Q_r(t,x)$.
We call $r$ the radius of $Q=Q_r(t,x)$.

We require a quite mild regularity
assumption on $a^{ij}$. They are assumed to be measurable in
$t$ and $x'$, and almost VMO with respect to $x''$. More
precisely, we impose the following assumption
in which $\theta>0$ will be specified later and
$R_{0}>0$ is a fixed number. Set
$$
\tr a=a^{11}+a^{22}.
$$

\begin{assumption}[$\theta$]
                                            \label{assump2}
For any
$Q=(s,t)\times B'\times B''\in\bQ$ with radius
$\rho\le R_0$
\begin{equation}
                                          \label{9.30.2}
 \dashint_Q|a (r,x)- a_{B''} (r,x')|\,dx\,dr
\le \theta,
\end{equation}
\begin{equation}
                                          \label{11.10.1}
\dashint_Q|\tr   a_{B''} (r,x')-
 \tr  a_{B'\times B''}(r) |\,dx\,dr
\le \theta,
\end{equation}
 where
$$
   a_{B''} (r,x')=
\dashint_{B''} a (r,x)\,dx'',
$$
$$
a_{B'\times B''}(r)= \dashint_{B'\times B''} a (r,x)\,dx
=\dashint_{B' }  a_{B''} (r,x')\,dx' .
$$
\end{assumption}

Observe that if $a(t,x)$ is independent of $x''$,
then the left-hand side of
 \eqref{9.30.2} is zero. It can be made as close to zero as we wish
on the account of $R_{0}$
if $a(t,x)$ is continuous with respect to $x''$
uniformly with respect to $(t,x')$. Also if
  $\tr a(t,x)$ depends only on $t,x''$ (for instance, constant),
then the left-hand side of 
 \eqref{11.10.1} is zero. It can be made as close to zero as we wish
on the account of $R_{0}$
if $\tr a(t,x)$ is continuous with respect to $x'$
uniformly with respect to $(t,x'')$.

In the following theorems we use the constants (small)
$\gamma_{0}=\gamma_{0}(\delta)\in(0,1/2)$ and $\theta=\theta(p,q,d,\delta)>0$
and (large)
$N=N(p,q,d,\delta,K,  R_{0})$
which will be determined later.

\begin{theorem}
                                           \label{theorem 11.10.1}
Let $p\in(2,2+\gamma_{0})$, $q>(pd-2p)/(2p-4) $  and let Assumption
\ref{assump2} ($\theta$) be  satisfied. Take $(t_{0},x_{0})\in\bR^{d+1}$.
Then
there exists a probability space $(\Omega,\cF,P)$
carrying a $d$-dimensional Wiener process $ (w_{t},\cF_{t}) $,
$t\geq0$,
and there exists a solution $x_{t}$ of the It\^o
equation
\begin{equation}
                                             \label{11.10.2}
x_{t}=x_{0}+\int_{0}^{t}a^{1/2}(t_{0}+s,x_{s})\,dw_{s}+
\int_{0}^{t}b(t_{0}+s,x_{s})\,ds,\quad t\geq0,
\end{equation}
such that it possesses the following
property  (a): for any nonnegative Borel $f(t,x)$ we have
\begin{equation}
                                             \label{11.10.3}
E\int_{0}^{\infty}e^{-t}f(t,x_{t})\,dt\leq N\|f\|_{L_{p,q}};
\end{equation}
(b) for any $u\in W^{1,2}_{p,q}$, $\lambda\in \bR$, 
and bounded $(\cF_{t})$-stopping time $\tau$ we have
\begin{equation}
                                             \label{11.10.4}
u(t_{0},x_{0})= 
E\int_{0}^{\tau}(\lambda u-Lu)(t_{0}+t,x_{t})
e^{-\lambda t-\phi_{t}}\,dt
+ E u(t_{0}+\tau,x_{\tau})e^{-\lambda  \tau -\phi_{\tau}},
\end{equation}
where
$$
\phi_{t}=\int_{0}^{t}c(t_{0}+s,x_{s})\,ds.
$$
\end{theorem}

\begin{definition}
                                    \label{definition 11.11.1}
We call any solution of  \eqref{11.10.2} 
(on any probability space) ``good''
if for some $p,q$ as in Theorem 
\ref{theorem 11.10.1} and
 any $T\in(0,\infty)$ there exists a constant $N$ such that
\begin{equation}
                                                      \label{27.1.1}
E\int_{0}^{T} f(s,x_{s})\,ds\leq N\|f\|_{L_{p,q}}
\end{equation}
  for any nonnegative Borel $f $.
\end{definition}

Observe that the solution existence of which is asserted
in Theorem 
\ref{theorem 11.10.1} is ``good''. Also note that
condition \eqref{27.1.1} {\em means\/} that there exists
a function $G(s,y)\geq0$ such that  
for any nonnegative Borel $f $ 
$$
E\int_{0}^{\infty} f(s,x_{s})\,ds=
\int_{0}^{\infty}\int_{\bR^{d}}G(s,y)f(s,y)\,dyds
$$
and for any $T\in(0,\infty)$ we have
\begin{equation}
                                  \label{1.27.4}
\int_{\bR^{d-2}}\Big( \int_{0}^{T}
\int_{\bR^{2}}|G|^{p'}\,dx'dt\big)^{q'/p'}\,dx''\leq N^{q'},
\end{equation}
where $N$ is taken from \eqref{27.1.1}, $p'=p/(p-1)$,
$q'=q/(q-1)$.

\begin{theorem}
                                       \label{theorem 11.11.1}
Under the assumptions of Theorem \ref{theorem 11.10.1} 
all ``good'' solutions of \eqref{11.10.2}  
 have the same finite-dimensional
distribution (weak uniqueness of solutions of
\eqref{11.10.2}).
\end{theorem}

The restrictions on $q$ come from the following.

\begin{theorem}
                                   \label{theorem 11.10.1}
Let   
one of the following conditions be satisfied:

(a) $d\geq 4$, $p>2$, $q>(pd-2p)/(2p-4)$,

(b) $d=3$, $2<p<4$, $q>p/(2p-4)$,

(c) $d=3$, $p\geq 4$, $q\in(1,\infty)$.

Then for any $u\in W^{1,2}_{p,q}$ and $(t,x)\in \bR^{d+1}$
we have
\begin{equation}
                                             \label{11.10.5}
|u(t,x)|\leq N\|I_{t}(\Delta u+\partial_{t}u-u)\|_{L_{p,q}}
\leq N\|u\|_{W^{1,2}_{p,q}},
\end{equation}
where $N=N(p,q,d)$ and $I_{t}(s,y)=I_{0}(s)=I_{(t,\infty)}(s)$.

\end{theorem}

First we need a lemma.
\begin{lemma}
                                    \label{lemma 1.30.1}
For $\lambda>0$ and fixed $\alpha\geq1$ introduce
$$
I(\lambda)=\int_{0}^{\infty}\frac{s}{(1+s^{2})^{\alpha/2}}
e^{-\lambda s}\,ds.
$$
Then $I(\lambda)$ is bounded for $\alpha>2$, and
 $\lambda^{\beta}I(\lambda)\to 0$ as $\lambda\downarrow 0$
for any $\beta>2-\alpha$ if $2\geq \alpha\geq 1$. 

\end{lemma}

Proof. The assertion in case $\alpha>2$ is obvious.
If $\alpha=1$, as is easy to see after the substitution
$s=t/\lambda$,
$\lambda I(\lambda)\to 1$ as $\lambda\downarrow 0$.
In the remaining case
  $2\geq \alpha>1$ and we integrate by parts to get
$$
I(\lambda)=-\frac{1}{\lambda}\int_{0}^{\infty}
(1-e^{-\lambda s})g(s)\,ds,\quad g(s)=\frac{d}{ds}
\frac{s}{(1+s^{2})^{\alpha/2}}.
$$
Clearly, $|g(s)|\leq N(1+s)^{-\alpha}$, where and below
by $N$ we denote constants depending
only on $\alpha$. Also, for any $\kappa\in(0,1]$, 
$1-e^{-\lambda s}\leq N\lambda^{\kappa}
s^{\kappa}$, which with $\kappa< \alpha-1$ allows us to write
$$
|I(\lambda)|\leq N\lambda^{\kappa-1}\int_{0}^{\infty}
(1+s)^{\kappa-\alpha}\,ds.
$$
The last integral is finite since $\kappa-\alpha<-1$,
$$
\nlimsup_{\lambda\downarrow0}\lambda^{1-\kappa}
|I(\lambda)|<\infty,
$$  
 and since for any $\beta$
such that $\beta>2-\alpha$ there is a $\kappa\in(0,\alpha-1)$
such that $\beta>1-\kappa$, the lemma is proved.

{\bf Proof of Theorem \ref{theorem 11.10.1}}. We may assume that $u\in C^{\infty}_{0}$
and $(t,x)=(0,0)$.
In that case set $f=u-\Delta u-\partial_{t}u$. Then
for a constant $c_{d}$
$$
u(0,0)=c_{d}\int_{\bR^{d-2}}\Big(
\int_{0}^{\infty}\int_{\bR^{2}}t^{-d/2}e^{-t-|x|^{2}/(4t)}
f(t,x)\,dx'dt\Big)\,dx''
$$
$$
\leq c_{d}\int_{\bR^{d-2}}g(x'')h(x'')\,dx'',
$$
where
$$
g(x'')=\Big(\int_{\bR^{3}}I_{0}(t )|f(t,x)|^{p}\,dx'dt\Big)^{1/p},
$$
$$
h(x'')=\Big(\int_{\bR^{3}}I_{0}(t )t^{-p'd/2}e^{-p't-p'|x|^{2}/(4t)}
\,dx'dt\Big)^{1/p'},
$$
and $p'=p/(p-1)$.
Note that for $a,b>0$
$$
 at+2b|x|^{2}/t\geq b|x|^{2}/t+2\sqrt{ab}|x|
\geq  b|x|^{2}/t+\sqrt{ab}|x' |+\sqrt{ab}|x''|,
$$
so that $e^{-p't-p'|x|^{2}/(4t)}\leq
e^{-\nu|x|^{2}/t-\mu|x |-\mu|x''|}$
for some $\mu,\nu>0$ depending only on $p$.

Taking this into account and then replacing $t$ in the integral
defining $h$ the variable $t$ with $ |x|^{2}s$ we find that
$$
h^{p'}(x'')\leq Ne^{-\mu |x''| }\int_{\bR^{2}}
|x|^{2-p'd}e^{-\mu |x' | }\,dx'
=Ne^{-\mu |x''| }|x''|^{ 4 -p' d }I(\mu|x''|),
$$
where $I$ is taken from Lemma \ref{lemma 1.30.1}
with $\alpha=p'd-2$ and the last equality is obtained
 by setting $x'=|x''|y$ and using polar coordinates.

Now we split the rest of the proof into three cases
according to our alternative assumptions.

(a) If $d\geq 4$ and $p>2$, then $1<p'<2$ and
$\alpha\in (d-2,2d-2)\subset(2,\infty)$. In that case, by Lemma \ref{lemma 1.30.1}
the function $I(\mu|x''|)$ is bounded and 
\begin{equation}
                                           \label{1.30.3}
h(x'')\leq Ne^{-\mu |x''|/p' }|x''|^{ 4/p' - d }.
\end{equation}

It follows that
$$
u(0,0)\leq N\|I_{0}f\|_{L_{p,q}}\Big(\int_{\bR^{d-2}}
e^{-q'p' |x''|/\mu }|x''|^{q'(4/p'- d)}\,dx''\Big)^{1/q'},
$$
where $q'=q/(q-1)$. The last integral is finite iff $p>2$
and $q>(pd-2p)/(2p-4)$ and this yields 
\eqref{11.10.5}. 

(b) If $d=3$ and $2<p<4$, then $4/3<p'<2$
and $\alpha=3p'-2\in
(2,4)$. In that case we finish the argument
as in (a).

(c) If $d=3$ and $p\geq 4$, then $1<p'\leq 4/3$ and
$\alpha=3p'-2\in(1,2]$. In that case
$|x''|^{ 4 -3p'   }I(\mu|x''|)=|x''|^{ 2 -\alpha }I(\mu|x''|)$
blows up at the origin slower than $|x|^{-\varepsilon}$
for any $\varepsilon>0$ and
$$
u(0,0)\leq N\|I_{0}f\|_{L_{p,q}}\Big(\int_{\bR^{d-2}}
e^{-q'p' |x''|/\mu }|x''|^{-\varepsilon}\,dx''\Big)^{1/q'}
$$
(with a different arbitrary small $\varepsilon>0$).
This proves the theorem. 

Recall that given $(t_{0},x_{0})$ a nonnegative function
$G(s,y)$ is called a Green's function of the operator
$L $ with pole at $(t_{0},x_{0})$ if the equality
\begin{equation}
                                  \label{1.27.2}
u(t_{0},x_{0})=-
\int_{0}^{\infty}\int_{\bR^{d}}G(s,y) L u(s,y)\,dyds
\end{equation}
holds for all $u\in C^{\infty}_{0}$.
We call $G$ a ``good'' Green's function if, 
for each $T\in(0,\infty)$, the left-hand side of
  \eqref{1.27.4} is finite
($p$ and $q$ are always taken from Theorem \ref{theorem 11.10.1}).
\begin{theorem}
                                        \label{theorem 1.27.1}
Let assumption of  Theorem
\ref{theorem 11.10.1} be satisfied 
and $(t_{0},x_{0})\in\bR^{d+1}$. Then there exists
a unique ``good'' Green's function of $L$
with pole at $(t_{0},x_{0})$.

\end{theorem}

The existence part in this theorem is just
a simple consequence of Theorem \ref{theorem 11.10.1}.
Indeed,  fix $T\in(0,\infty)$. Then
$$
E\int_{0}^{T}e^{-\phi_{t}}f(t,x_{t})\,dt
$$
is a nonnegative linear bounded functional on $L_{p,q}$.
Hence there exists $G\geq0$ for which
the left-hand side of
  \eqref{1.27.4} is finite and
$$
E\int_{0}^{T}e^{-\phi_{t}}f(t,x_{t})\,dt=
\int_{0}^{T}\int_{\bR^{d}}G(s,y)f(s,y)\,dyds.
$$
Obviously, $G$ is independent of $T$ and
\begin{equation}
                                  \label{1.27.6}
E\int_{0}^{\infty}e^{-\phi_{t}}f(t,x_{t})\,dt=
\int_{0}^{\infty}\int_{\bR^{d}}G(s,y)f(s,y)\,dyds
\end{equation}
for any nonnegative Borel $f$ or for any $f$
for which at least one side is finite.
By taking $f=-Lu$, where $u\in C^{\infty}_{0}$
and using \eqref{11.10.4} with $\tau$ so large that
$u(t_{0}+\tau+s,x)=0$ for $s\geq0$, we immediately get
\eqref{1.27.2}.

 The most important ingredient in the proof
of Theorems \ref{theorem 11.10.1},
 \ref{theorem 11.11.1}, and the uniqueness part in
Theorem 
\ref{theorem 1.27.1} is the following
result which we prove in Section \ref{sec4}.
It generalizes one of the main results of \cite{DK_10}
in the respect that the continuity of $\tr a$
with respect to $x''$ uniform with respect to $(t,x'')$ is replaced
by a kind of VMO condition.
Recall that $\gamma_{0}=\gamma_{0}(\delta)\in(0,1/2)$,
$\theta=\theta(p,q,d,\delta)>0$
and 
$N=N(p,q,d,\delta,K,  R_{0})$ are mentioned before Theorem
\ref{theorem 11.10.1}.
\begin{theorem}
                                      \label{theorem 11.3.1}
 
There exists    $\lambda_{0}>0$,
depending only on
$p,q,d,\delta,K,  R_{0}$,
such that if $p\in(2,2+\gamma_{0})$, $q\in(2,\infty)$,
and Assumption
\ref{assump2} ($\theta$) is satisfied, then

(i)
For any $u\in W^{1,2}_{p,q} $  and  
$\lambda\geq \lambda_{0}$
\begin{equation}
                                                   \label{11.3.2}
\lambda\|u\|_{L_{p,q}}+
\sqrt{\lambda}\|Du\|_{L_{p,q}}+
\|D^{2}u,\partial_{t}u\|_{L_{p,q}}
\leq N\|\lambda u-Lu\|_{L_{p,q}};
\end{equation}

(ii) For any $\lambda\geq \lambda_0$ and $f\in L_{p,q} $, there exists a unique solution
$u\in W^{1,2}_{p,q} $ of equation \eqref{parabolic} in
$\bR^{d+1} $.
 
\end{theorem}
\begin{remark}
                                             \label{remark 1.27.3}

 Theorem \ref{theorem 11.3.1} and formula \eqref{11.10.4}
 imply the solvability of the Cauchy
problem as in \cite{Kr_07}.  Indeed, if we take $f(t,x)=0$
for $t\geq T$, then for the solution $u$ of \eqref{parabolic}
we get  from  \eqref{11.10.4} 
that 
$|u(t_{0},x_{0})|\leq e^{-\lambda(S-t_{0})}\sup|u|$
for any $t_{0}\geq T$ and $S>T$,
where the last sup is finite due to Theorem \ref{theorem 11.10.1}.
By sending $S\to \infty$ we conclude that $u(t,x)=0$
for $t\geq T$ and hence $u$ is a solution of the Cauchy problem
for \eqref{parabolic} for $t<T$ with terminal data
$u=0$.

\end{remark} 

{\bf Proof of Theorem \ref{theorem 11.10.1}}.
 First note that there is a sequence of
operators $L_{n}$ with infinitely differentiable coefficients
$a_{n},b_{n},c_{n}$
which satisfy the same assumptions as the original $a,b,c$
and which converge to $a,b,c$
almost everywhere in $\bR^{d+1}$. To see that it suffices to use mollifiers
with nonnegative kernels. Next, take a probability
space $(\Omega',\cF',P')$ carrying a $d$-dimensional Wiener
process $w'_{t}$, $t\geq0$, and define $x^{n}_{t}$, $t\geq0$, as  
unique solutions of
$$
x^{n}_{t}=x_{0}+\int_{0}^{t}a^{1/2}_{n}(t_{0}+s,x^{n}_{s})\,dw'_{s}+
\int_{0}^{t}b_{n}(t_{0}+s,x^{n}_{s})\,ds,\quad t\geq0.
$$

We know (consequence of the parabolic Aleksandrov
estimates and Skoro\-khod's embedding method,
see, for instance, Section 2.6 of \cite{Kr_77})
that the collection of the distributions
of $x^{n}_{\cdot}$ on $C([0,\infty),\bR^{d})$
is tight and for any weakly convergent subsequence
there exist probability space $(\Omega,\cF,P)$
carrying a $d$-dimensional Wiener process $ (w_{t},\cF_{t}) $,
$t\geq0$,
and   a solution $x_{t}$ of  
equation \eqref{11.10.2} such that a subsequence
of distributions of $x^{n}_{\cdot}$ converges to the distribution
of $x_{\cdot}$. Furthermore, along the same
subsequence, for any Borel $\bR^{m}$-valued bounded
function $g(t,x)$, the distributions on 
$$
C([0,\infty),\bR^{d})\times C([0,\infty),\bR^{m})
$$
of 
 \begin{equation}
                                               \label{11.11.1}
( x^{n}_{\cdot},\int_{0}^{\cdot}g(s,x^{n}_{s})\,ds)
\end{equation}
weakly converge to
the distribution of 
$$
( x _{\cdot},\int_{0}^{\cdot}g(s,x_{s})\,ds).
$$
 Our first goal is to show
that $x_{\cdot}$ possesses   property (a)  with $p,q$ as in 
the statement of the theorem for which it holds 
automatically that $q>p$ because $\gamma_{0}\in(0,1/2)$

Note that for any $f\in C^{\infty}_{0}$ 
there exists a smooth solution $u^{n}$
of $\lambda_{0} u^{n}-L_{n}u^{n}=f$ for which
Theorem \ref{theorem 11.3.1} is valid.
In case $c_{n}\equiv0$, by It\^o's formula
$$
u^{n}(t_{0},x_{0})=E\int_{0}^{\infty}e^{-t}
\big[f-(\lambda_{0}-1)u^{n}
\big](t_{0}+t,x^{n}_{t})\,dt,
$$
which implies by Theorem \ref{theorem 11.3.1}
and Theorem \ref{theorem 11.10.1} that
$$
E\int_{0}^{\infty}e^{-t}f(t_{0}+t,x^{n}_{t})\,dt
\leq (1+\lambda_{0})\sup_{\bR^{d+1}}| u^{n}|\leq N\|f\|_{L_{p,q}}.
$$
By passing to the limit we see that \eqref{11.10.3}
holds 
for $f\in C^{\infty}_{0}$ and then a standard argument
shows that it also holds for all Borel nonnegative $f$.

Next, take $u\in W^{1,2}_{p,q}$ and
let $u^{(\varepsilon)}$ be mollified $u$ with smooth kernels
supported in $(-\varepsilon,\varepsilon)\times B_{\varepsilon}$
such that $u^{(\varepsilon)}\to u$ in $W^{1,2}_{p,q}$
as $\varepsilon\downarrow0$.
In particular, $u^{(\varepsilon)}\to u$ uniformly and
$\lambda  u^{(\varepsilon)} -L u^{(\varepsilon)}\to 
\lambda  u  -L u$
in $L_{p,q}$. Then after writing It\^o's formula
for $u^{(\varepsilon)}(t_{0}+t,x_{t})\exp[ -\lambda 
 t-\phi_{t}]$,
passing to the limit, and using \eqref{11.10.3}
(recall that $\tau$ is bounded) we 
immediately arrive at \eqref{11.10.4}.
The theorem is proved.

{\bf Proof of Theorem \ref{theorem 11.11.1}}.
 The last paragraph in the above proof convinces us that
\eqref{11.10.4} holds for any ``good'' solution of \eqref{11.10.2}.
Moreover,   for   $f\in C^{\infty}_{0}$
by Theorem \ref{theorem 11.3.1}
there is a unique solution $u\in W^{1,2}_{p,q}$
of $\lambda_{0} u -L u^{n}=f$. Since
$$
\lambda_{0} (u^{n}-u)-L^{n}(u^{n}-u)=(L-L_{n})u,
$$
estimate \eqref{11.3.2} and the dominated
convergence theorem
imply that $u^{n}\to u$
in $W^{1,2}_{p,q}$, and then $u^{n}\to u$ uniformly
on $\bR^{d+1}$. Furthermore, it is a classical fact that
if $T\in(0,\infty)$ and $f(t_{0}+t,x)=0$ for $t\geq T$, then  
$u^{n}(t_{0}+t,x)=0$ and hence $u(t_{0}+t,x)=0$
for $t\geq T$. By using \eqref{11.10.4} with $\tau=T$
we conclude that
$$
E\int_{0}^{T}f(t_{0}+t,x_{t})e^{-\lambda_{0}t-\phi_{t}}
\,dt=u(t_{0},x_{0})
$$
and hence the left-hand side is independent of the choice
of solution of \eqref{11.10.2} provided $f\in C_{0}^{\infty}$
and $f(t_{0}+t,x)=0$ for $t\geq T$.
By usual measure theoretic argument one shows that
the independence holds for any Borel $f$ which is bounded
on $[0,T]\times \bR^{d}$. In particular,
$$
E\int_{0}^{T}c(t_{0}+t,x _{t})
e^{ -\phi _{t}}\,dt=1-E\exp\Big(-\int_{0}^{T}
c(t_{0}+t,x_{t})\,dt\Big)
$$
is independent of the choice
of solution of \eqref{11.10.2}
for any bounded $c$.
This easily implies weak uniqueness. The theorem
is proved.

\begin{corollary}
                                     \label{corollary 11.11.1}
The property of weak uniqueness and the tightness of distributions
of $x^{n}_{\cdot}$ obviously imply that the whole
sequence of distributions of $x^{n}_{\cdot}$
weakly converges to the distribution of $x_{\cdot}$.
\end{corollary}

We fix $\theta$ from Theorem \ref{theorem 11.10.1}
and suppose that Assumption
\ref{assump2} ($\theta$) is satisfied.
Next, denote by $\Omega $ the set
of $\bR^{d+1}$-valued functions
$\omega=\omega_{s}=(t+s,x_{s}),s\in[0,\infty)$,
such that $x_{\cdot}\in 
C([0,\infty),\bR^{d})$ and, if $\omega
=\{(t+\cdot,x_{\cdot})$ and $s\in[0,\infty)$, we set  $t_{s}(\omega)=t+s$ and $x_{s}(\omega)=x_{s}$.
As usual the argument $\omega$ is almost always dropped.
Introduce $\cN_{r}=\sigma\{(t_{s},x_{s}),0\leq s\leq r\}$,
$\cN_{\infty}=\sigma\{(t_{s},x_{s}),0\leq s<\infty\}$.

If $x_{t}$, $t\geq0$, is a solution of
\eqref{11.10.2} on a probability space, then
the function $(t_{0}+t,x_{t})$, $t\geq0$,
is an $\Omega$-valued, $N_{\infty}$-measurable
random  variable. If the solution is ``good'',
its distribution on $\Omega$
we denote by   $P_{t_{0},x_{0}}$. Obviously $P_{t_{0},x_{0}}$
are defined for any $(t_{0},x_{0})$.

\begin{theorem}
                                          \label{theorem 11.13.1}
The triplet consisting of 
$\Omega$, the family $\cN_{t}$, $t\geq0$,
and the family $P_{t ,x }$, $(t,x)\in\bR^{d+1}$,
is a strong Markov process which is strong Feller
in the sense that for any $T\in(0,\infty)$
and Borel bounded $f(x)$ on $\bR^{d}$ the function
\begin{equation}
                                               \label{11.13.1}
v(t,x)=\int_{\Omega}f(x_{T-t})\,P_{t,x}(d\omega)
\end{equation}
is $\alpha$ H\"older continuous in $(t,x)$ such that
$t<T$, where $\alpha\in(0,1)$ depends only on $d$ and 
$\delta$.

\end{theorem}

Proof. Take $f\in C^{\infty}_{0}$ and set
$$
u(t,x)=\int_{\Omega}\int_{0}^{\infty}f(t,x_{t})
e^{-\lambda_{0}t}\,dt\,P_{t,x}(d\omega).
$$
Then let $v$ be   any $W^{1,2}_{p,q}$-solution
of $\lambda_{0} v-Lv=f$ with $c\equiv0$ (and $p,q$
as in Theorem \ref{theorem 11.10.1}).
By \eqref{11.10.4}, in which we take
$\lambda=\lambda_{0}$, $\tau=T$ and let $T
\to \infty$, we see that
$$
v(t_{0},x_{0})=E\int_{0}^{\infty}
f( t_{0}+t,x_{t})
e^{-\lambda_{0}t}\,dt.
$$
It follows   that $u=v$ at $(t_{0},x_{0})$
and at every other point as well.
Now the strong Markov property follows directly from
\eqref{11.10.4}.

The strong Feller property for the Markov
process generated by  smooth coefficients $a_{n},b_{n}$, taken from the proof of Theorem \ref{theorem 11.10.1},
is a classical result. The fact that the H\"older
exponent and constants are under control
independent
of the smoothness of $a_{n},b_{n}$ is the Krylov-Safonov
result. Furthermore, Corollary \ref{corollary 11.11.1}
says that the functions \eqref{11.13.1}
corresponding to $a_{n},b_{n}$ converge to
$v$ if $f$ is continuous. However,
the above mentioned estimates of the  H\"older continuity
do not involve anything from $f$ apart from 
 $\sup|f|$ and this and a usual
measure theoretic argument proves
the  H\"older continuity for any  Borel
bounded $f$. The theorem is proved.

{\bf Proof of Theorem \ref{theorem 1.27.1}}.
We only need to prove uniqueness assuming without loss of
generality that $(t_{0},x_{0})=0$. Take $T\in(0,\infty)$
and observe that, owing to obvious approximations
and the assumption concerning \eqref{1.27.4},
\eqref{1.27.2} holds with $(t_{0},x_{0})=0$
 for ``good'' Green's functions
not only for $u\in C^{\infty}_{0}$ but also for
$u\in W^{1,2}_{p,q}$ if $u(t,x)=0$ for $t\geq T$.

Then take
$f\in L_{p,q}$ such that $f(t,x)=0$ if $t\geq T$,
let $u\in W^{1,2}_{p,q}$ be a unique solution
of \eqref{parabolic} with $\lambda=\lambda_{0}$,
and set $v(t,x)=e^{\lambda_{0}t}u(t,x)$. Observe that
$v\in W^{1,2}_{p,q}$ since $u(t,x)=0$ for $t\geq T$
(see Remark \ref{remark 1.27.3}) and
$
Lv=-e^{ \lambda_{0}t}f
$.
  Then  
$$
 \int_{0}^{T}\int_{\bR^{d}}G(s,y)
e^{ \lambda_{0}s}f(s,y)\,dyds
$$
 is the same for all ``good'' Green's functions,
because it is equal to $v(0)$ by the above. The arbitrariness
of $f\in L_{p,q}$ and the assumption that the left-hand side of 
\eqref{1.27.4} is finite bring the proof to an end.

\section{Preliminary results}
                                            \label{sec2}
We first consider equations in $\bR\times \bR^2$ with measurable
coefficients. 
\begin{lemma}
                                \label{lem2.1}
Let   $d=2$ and
$$
Lu=\partial_{t} u +\sum_{i,j=1}^2 a^{ij}(t,x)D_{ij}u. 
$$
Assume that $\tr a$ depends only on $t$. Then there exists
 a  $\gamma_0=\gamma_0(\delta)>0$ such that for
any
$p\in (2-\gamma_0,2+\gamma_0)$, $u\in W^{1,2}_p(\bR^3 )$, and
$\lambda\ge 0$, we have
$$
\|D^{2}u\|_{L_{p}(\bR^3 )}+\|\partial_{t} u\|_{L_{p}(\bR^3 )}
+\sqrt{\lambda}
\|Du\|_{L_{p}(\bR^3 )}
$$
\begin{equation}
                                            \label{2.52}
+\lambda\|u\|_{L_{p}(\bR^3 )}
\leq N\|Lu-\lambda u\|_{L_{p}(\bR^3 )},
\end{equation} where $N=N( \delta,p)$. Moreover for any $\lambda>0$ and
$f\in L_p(\bR^3 )$ there exists a unique $u\in W^{1,2}_p(\bR^3 )$
solving $Lu-\lambda u=f$  in $\bR^3 $.
\end{lemma}

This lemma is proved in \cite{DK_10} as Lemma 3.1.

\begin{lemma}
                                               \label{lemma 10.28.1}
For any $ p,  q\in(1,\infty)$ and $\varepsilon>0$
there exist $N_{0}=N_{0}(d,p,q)$ and
 $N=N(d,p,q,\varepsilon)$ such that
for any $u$ such that
$u,Du,D^{2}u\in L_{p,q}$ and $\lambda\geq0$ we have
\begin{equation}
                                                  \label{11.4.2}
\lambda\|D^{2}u\|_{L_{p,q}}+\lambda^{1/2}\|Du\|_{L_{p,q}}
+\|D^{2}u\|_{L_{p,q}}
\leq N \|\lambda u-\Delta u\|_{L_{p,q}},
\end{equation}
\begin{equation}
                                                  \label{11.4.3}
\|D_{x'x''}u\|_{L_{p,q}}
\leq \varepsilon \|D_{x'x'}u\|_{L_{p,q}}+
N \|D_{x''x''}u\|_{L_{p,q}}.
\end{equation}

 \end{lemma}

Proof. First observe that if \eqref{11.4.2} is true, then
$$
\|D_{x'x''}u\|_{L_{p,q}}
\leq N\|D_{x'x'}u\|_{L_{p,q}}+
N \|D_{x''x''}u\|_{L_{p,q}}
$$
and \eqref{11.4.2} follows owing to the different homogeneity
of the above terms with respect to 
scalings in $x'$.

Owing to the possibility of mollification, while proving
\eqref{11.4.2} we may assume that $\partial_{t}u,\partial^{2}_{t}u
\in L_{p,q}$.
Then the possibility to use scalings in $t$ shows
that to prove \eqref{11.4.2} it suffices to show that
\begin{equation}
                                                  \label{11.4.4}
\lambda\|D^{2}u\|_{L_{p,q}}+\lambda^{1/2}\|Du\|_{L_{p,q}}
+\|D^{2}u\|_{L_{p,q}}
\leq N \|\lambda u-\Delta u-\partial^{2}_{t}u\|_{L_{p,q}}.
\end{equation}

That \eqref{11.4.4} holds for sufficiently large
$\lambda$ (with $N$ independent of $\lambda$) follows
from Theorem 5.5 of \cite{DK_18}
as a very particular case. Then the fact that it holds for 
any $\lambda\geq0$ follows by scaling in $(t,x)$.
The lemma is proved.

An immediate corollary of Lemmas \ref{lem2.1} 
and \ref{lemma 10.28.1}
is the following estimate.

\begin{corollary}
                                            \label{corollary 2.2}
Assume that $\tr a$
depends only on $(t,x'')$.
Let    
$$
Lu= u_t+\sum_{i,j=1}^d a^{ij}( t, x)D_{ij}u.
$$
  Then for any
$p\in (2-\gamma_0,2+\gamma_0)$, where $\gamma_0$ is taken from
 Lemma \ref{lem2.1},  any
$q\in(1,\infty)$, and any $u\in W^{1,2}_{p,q} $ 
and $\lambda\ge
0$, we have
$$
\lambda\|u\|_{L_{p,q} }+\sqrt{\lambda}
\|Du\|_{L_{p,q} }+\|D^{2}u\|_{L_{p,q} }
+\|\partial_{t} u\|_{L_{p,q} }
$$
\begin{equation}
                                            \label{2.57}
\leq N\|Lu-\lambda u\|_{L_{p,q} }
+N\|D^{2}_{x''}u\|_{L_{p,q} },
\end{equation}
where $N=N(\delta,d,p,q )$.
\end{corollary}

Proof.
We first fix $x''$ and apply Lemma \ref{lem2.1} to get
$$
\lambda^{q}\|u(\cdot,\cdot,x'')\|^q_{L_{p}(\bR^3 )}
+\|D_{x'}^{2}u(\cdot,\cdot,x'')\|^q_{L_{p}(\bR^3 )}
+\|\partial_{t}u(\cdot,\cdot,x'')\|^q_{L_{p}(\bR^3 )}
$$
\begin{equation}
                                            \label{3.50}
\leq N\|  \sum_{i,j=1}^2 a^{ij}   D_{ij}u(\cdot,\cdot,x'')
+\partial_{t}u  (\cdot,\cdot,x'') -\lambda
u(\cdot,\cdot,x'')\|^q_{L_{p}(\bR^3 )}.
\end{equation}
Upon integrating \eqref{3.50} with respect to $x''$ we arrive at
$$
\lambda\|u\|_{L_{p,q} }+\|D_{x'}^{2}u\|_{L_{p,q} }
+\|\partial_{t}u\|_{L_{p,q} }
$$
\begin{equation}
                                            \label{3.55}
\leq N\|Lu-\lambda u\|_{L_{p,q} }
+\|D_{xx''}u\|_{L_{p,q} }.
\end{equation}

  By using Lemma \ref{lemma 10.28.1}, we  deduce from
\eqref{3.55} that
$$
\lambda\|u\|_{L_{p,q} }
+\|D^{2}u\|_{L_{p,q} }
+\|\partial_{t}u\|_{L_{p,q} }
$$
\begin{equation*}
\leq N\|Lu-\lambda u\|_{L_{p,q} }
+\|D_{x''x''}u\|_{L_{p,q} }.
\end{equation*}
To estimate  $\|Du\|_{L_{p,q} }$, we 
again use Lemma \ref{lemma 10.28.1} showing that 
$$
\sqrt \lambda \|Du\|_{L_{p,q}}\le N \lambda \|u\|_{L_{p,q}}
+N\|D^2u\|_{L_{p,q}}.
$$
The corollary is proved.

\section{Proof of Theorem 
\protect\ref{theorem 11.3.1}}
                                            \label{sec4}
  We consider the
operator
$$
L u= \partial_{t}u+a^{ij}D_{ij}u+b^{i}D_iu-cu,
$$
where $(a^{ij})$ satisfy Assumption \ref{assump2} ($\theta$) with
some $\theta>0$ to be specified later. First we deal with 
assertion (i) of Theorem \ref{theorem 11.3.1}.

Introduce
$$
\Var_{ Q}f=\dashint_{Q}|f-f_{Q}|^{2}\,dxdt.
$$
Take a nonnegative $\zeta\in C^{\infty}_{0}(Q_{R_{0}})$
whose square integrates to one.
\begin{lemma}
                                       \label{lemma 10.30.1}
For any $Q\in \bQ$ and $u$ such that $D_{x''}^{2}u\in L_{2}(Q)$,
we have
\begin{equation}
                                         \label{10.30.1}
\int_{\bR^{d+1}}\Var_{ Q}[D_{x''}^{2}(u\zeta(\cdot-(s,y))]\,dyds
\geq\Var_{ Q}v,
\end{equation}
where 
\begin{equation}
                                             \label{10.31.1}
v=(|D_{x''}^{2}u|^{2}+\chi|u|^{2})^{1/2},
\end{equation}
and $\chi=\chi(d,R_{0})$.
\end{lemma}

Proof. We know that
$$
2\Var_{ Q}f=\dashint_{Q}\dashint_{Q}|f(t_{1},x_{1})
 -f(t_{2},x_{2})|^{2}\,dx_{1}dx_{2}dt_{1}dt_{2}.
$$
We also know that
$$
\int_{\bR^{d+1}}|f(s,y)-g(s,y)|^{2}\,dyds
$$
$$
\geq\Big|\Big (\int_{\bR^{d+1}}|f(s,y) |^{2}\,dyds\Big)^{1/2}
-\Big (\int_{\bR^{d+1}}|g(s,y) |^{2}\,dyds\Big)^{1/2}\Big|^{2}.
 $$
Furthermore,
$$
\int_{\bR^{d+1}}|D_{x''}^{2}\big(u(t_{1},x_{1})\zeta
(t_{1}-s,x_{1}-y)\big)\big |^{2}\,dyds
$$
$$
=|D_{x''}^{2} u(t_{1},x_{1})|^{2}+\chi|u(t_{1},x_{1})|^{2},
$$
where 
$$
\chi=\int_{\bR^{d+1}}|D_{x''}^{2}\zeta|^{2}\,dyds.
$$
Upon combining the above we come to \eqref{10.30.1}.
The lemma is proved.

Then, we extend Theorem 5.1 of \cite{DK_10}
in which $u$ was required to have support
in a translate of $Q_{R_{0}}$. This assumption
in \cite{DK_10} was harmless because we could use
partitions of unity for the particular operators
under consideration. Here we do not know how to use
partitions of unity in the mixed norm setting. 

\begin{theorem}
                                            \label{thm4.1}
Let $\alpha=\alpha(d,\delta)\in(0,1)$ be the constant in
Theorem 4.4 of \cite{DK_10},   $\tau,\sigma \in (1,\infty)$,
$1/\tau+1/\sigma=1$. Take a $u\in W^{1,2}_{2,\loc}$,
introduce $v$ by \eqref{10.31.1}, and set $f=L  u$. Then under
Assumption
\ref{assump2} ($\theta$) there exists a positive constant 
$N$ depending only on
$d$,
$\delta$, and $\tau$ such that, for any
$(t_0,x_0)\in \bR^{d+1}$,
$r\in (0,\infty)$, and $\kappa\ge 4$,
$$
\Var_{ Q_r(t_0,x_0)}v\le N\kappa^{d+2}
\left(|f|^2+|Du |^2+|u|^2\right)_{Q_{\kappa r}(t_0,x_0)}
$$
\begin{equation}          
                             \label{eq13.5.05}
 +N\kappa^{d+2}\theta^{1/\sigma}
\left(|D^{2}u|^{2\tau}\right)_{Q_{\kappa r}(t_0,x_0)}^{1/\tau}+N\kappa^{-2\alpha}
 (|D^{2}_{x''}u |^2 )_{Q_{\kappa r}(t_0,x_0)}.
\end{equation}
 \end{theorem}

Proof. Without losing generality we assume that
$(t_0,x_0)=(0,0)$, then we
 fix   $\kappa\ge 4$  and $r\in
(0,\infty)$. 
Owing to the fact that $|Du|$ and $|u|$ enter the right-hand side
of \eqref{eq13.5.05}, we may also assume that $b\equiv 0$
and $c\equiv $0.

Next,
for $(s,y)\in\bR^{d+1}$ set
  $Q=Q[s,y]=(s_{1},s_{2})\times B'\times B''$ to be  $Q_{\kappa
r} $ if $\kappa r< R_0 $ and   $Q_{ R_0 }(s,y)$ if $\kappa
r\ge  R_0 $. For such $Q[s,y]$ we denote
$B'=B'[s,y], B''=B''[s,y]$. Recall the definitions
  given in Assumption \ref{assump2} and
set  
$$
a_{[s,y]}(t)=a_{ B'[s,y]\times  B''[s,y]}(t )
$$
$$
\sfa _{[s,y]}=\frac{ a _{B''[s,y]}}{\tr  a_{B''[s,y]}}
 \tr a_{[s,y]} ,\quad
\hat{f}_{[s,y]}= \sfa ^{ij }_{[s,y]}D_{ij }u+\partial_{t} u.
$$
Obviously, $\sfa$ depends only on $(t,x')$,
$\tr\sfa=\tr a_{ [s,y]}$
 depends only on $t$ and takes values
  between $2\delta$ and $2\delta^{-1}$  and for
  any $(s,y)$,
$$
\dashint_{Q[s,y]}|a - \sfa  _{[s,y]}|\,dx dt\leq N
\dashint_{Q[s,y]}| a\, \tr  a_{B''[s,y]}
- \sfa  _{[s,y]}\tr  a_{B''[s,y]} |\,dx dt
$$
\begin{equation}
                                                   \label{9.16.1}
\leq N\dashint_{Q[s,y]}|a - a _{B''[s,y]}|\,dx\,dt+
N\dashint_{Q[s,y]}|\tr  a_{B''[s,y]}-  \tr a_{ [s,y]} 
|\,dx\,dt\leq N\theta.
\end{equation}
 
Note that
$$
\hat{f}_{[s,y]} =  \left(\sfa^{ij }_{[s,y]} -
 a^{ij }\right) D_{ij }u + f,
$$
and for any values of the parameters
 $(s,y)$ and 
$$
w_{[s,y]}(t,x):=u(t,x)\zeta(t-s,x-y)
$$
we have
$$
\sfa^{ij}_{[s,y]}D_{x^{i}x^{j}}w_{[s,y]}+\partial_{t}w_{[s,y]} 
=\hat f_{[s,y]}\zeta(\cdot-(s,y))+2\sfa^{ij}_{[s,y]}D_{x^{i}}u
D_{x^{j}}\zeta(\cdot-(s,y))
$$
$$
+u\sfa^{ij}_{[s,y]}D_{x^{i}x^{j}}
\zeta(\cdot-(s,y))+u\partial_{t}\zeta(\cdot-(s,y))=:\check f_{[s,y]}.
$$

Since $\tr \sfa$ depends only on $t$, by Theorem 4.5 
of \cite{DK_10}
with an appropriate translation  
\begin{equation}							
                                \label{13.5.42}
\Var_{Q_r } D^{2}_{x''}w_{[s,y]}   \le N\kappa^{d+2}
\left(| \check f_{[s,y]}|^2\right)_{Q_{\kappa r}}
+N\kappa^{-2\alpha}
\left(|D^{2}_{x''}w_{[s,y]}|^2\right)_{Q_{\kappa r} },
\end{equation}
where $N$ and $\alpha$ depend only on $d$ and $\delta$.
By the definition of $\check f_{[s,y]}$,
\begin{equation}
                                          							\label{13.5.52}
\dashint_{Q_{\kappa r}} |\check f_{[s,y]} |^2 \, dx dt
\le N \dashint_{Q_{\kappa r} } |f|^2\zeta^{2}(\cdot-(s,y)) \, dx dt
+N I_{[s,y]}+NJ_{[s,y]},
\end{equation}
where
$$
I_{[s,y]} =
\dashint_{Q_{\kappa r} }
 |  \sfa _{[s,y]} - a |^{2}\,| D^{2}u|^{2} \zeta^{2}(\cdot-(s,y)) \, dx dt,
$$
$$
J_{[s,y]}=\dashint_{Q_{\kappa r} }
(|D\zeta |^{2}+|D^{2}\zeta |^{2}+|\partial_{t}\zeta|^{2}
 )(\cdot-(s,y))(|Du|^{2}+|u|^{2})\,dxdt.
$$
Observe that owing to the facts that $a$ and $\zeta$
are bounded functions and $\zeta$ is square integrable,
 and taking into account  H\"older's inequality, we have
$$
\int_{\bR^{d+1}}I_{[s,y]}\,dyds
=\dashint_{Q_{\kappa r} }\Big(\int_{\bR^{d+1}}
 | \sfa _{[s,y]} - a |^{2}
\zeta^{2}(\cdot-(s,y))\,dyds\Big)\,| D^{2}u|^{2}   \, dx dt
$$
\begin{equation}	
                                          						\label{10.31.3}
\leq NI_{1}^{1/\sigma}I_{2}^{1/\tau},
\end{equation}
where
$$
I_{1}=\dashint_{Q_{\kappa r}}\int_{\bR^{d+1}}
|  \sfa_{[s,y]} - a | 
\zeta^{2}(\cdot-(s,y))\,dyds\, dx dt,
$$
$$
I_2 = \dashint_{Q_{\kappa r} } |D^{2}u|^{2\tau} \, dx dt.
$$

Note that,
if $\kappa r<R_{0}$, 
we have $Q[s,y]=Q_{\kappa r}$ and
$\sfa_{[s,y]}$ is independent of $(s,y)$
and is constructed on the basis of $Q_{\kappa r}$.
In light of this and \eqref{9.30.2} in that case
$I_{1}\leq N\theta$.

However, if $\kappa r\geq R_{0}$, we know that $\zeta$
is supported in $Q_{R_{0}}$ and therefore,
\begin{equation}	
                                          						\label{10.31.4}
I_{1}=N(\kappa r)^{-d-2}\int_{\bR^{d+1}}\Big(
\int_{Q_{\kappa r}\cap Q_{R_{0}}(s,y)}|  \sfa _{[s,y]} - a | 
\zeta^{2}(\cdot-(s,y))\, dx dt\Big)\,dyds.
\end{equation}
Here the intersection is nonempty only if $|y|\leq \kappa r+2R_{0}$
and $-R_{0}^{2}\leq s\leq (\kappa r)^{2} $.
Such couples are occupying the volume less than $N(\kappa r)^{d+2}$
(since $\kappa r\geq R_{0}$). For any of those couples the interior integral in \eqref{10.31.4} is less than the integral over 
$Q_{R_{0}}(s,y)$ and is dominated by the
$\max \zeta^{2}$ times the volume
of $Q_{R_{0}}$ times the first
expression in \eqref{9.16.1} with $Q[s,y]=Q_{R_{0}}(s,y)$. It follows that
in this case $I_{1}\leq N\theta$
again.

Estimating the integrals with respect to $(s,y)$ of $J_{[s,y]}$
and of the last term
in \eqref{13.5.42}
is straightforward.
This together with  \eqref{13.5.42}-\eqref{10.31.3} yields
\eqref{eq13.5.05}. The theorem is proved.
 
Next we extract some consequences from Theorem
\ref{thm4.1} in terms of maximal and sharp functions
in parabolic setting. Recall that the maximal function of $u$
is defined by
$$
\bM u(t,x)=\bM |u|(t,x)=\sup_{\substack{Q\in\bQ\\
Q\ni(t,x)}}\dashint_{Q}|u|\,dyds.
$$
Obviously, at each point of $Q_{r}(t_{0},x_{0})$
and even of $Q_{\kappa r}(t_{0},x_{0})$ 
the right-hand side of \eqref{eq13.5.05}
is less than the expression which you get from it
by replacing all averages with the corresponding maximal functions.
On the other hand, the sharp function of $u$ is defined by 
$$
 u^{\#}(t,x) =\sup_{\substack{Q\in\bQ\\
Q\ni(t,x)}}\dashint_{Q}|u-u_{Q}|\,dyds.
$$
In addition
$$
(\Var_{ Q_r(t_0,x_0)}v)^{1/2}
\geq\dashint_{ Q_r(t_0,x_0)}|v-v_{ Q_r(t_0,x_0)}|\,dxdt.
$$
It follows from the above that 
owing to \eqref{eq13.5.05}   for any $Q\in \bQ$
and  $(t,x)\in Q$ we have
$$
\Big(\dashint_{Q}|u-u_{Q}|\,dyds\Big)^{2}
\leq 
N\kappa^{d+2}
\bM\left(|f|^2+|Du |^2+|u|^2\right) (t,x)
$$
$$
 +N\kappa^{d+2}\theta^{1/\sigma}
\left(\bM(|D^{2}u|^{2\tau})\right) ^{1/\tau}(t,x)+N\kappa^{-2\alpha}
 \bM(|D^{2}_{x''}u |^2 )(t,x) .
$$
 Obviously we can replace the left-hand side here with
$(v^{\#}(t,x))^{2}$ and get
 $$
(v^{\#} )^{2}
\leq 
N\kappa^{d+2}
\bM\left(|f|^2+|Du |^2+|u|^2\right) 
$$
\begin{equation}          
                             \label{eq13.5.05}
 +N\kappa^{d+2}\theta^{1/\sigma}
\left(\bM(|D^{2}u|^{2\tau})\right) ^{1/\tau} +N\kappa^{-2\alpha}
 \bM(|D^{2}_{x''}u |^2 ) .
\end{equation}
 
Now
we  remind the reader some well-known properties
of the $A_{p}$-weights ($p\in(1,\infty))$.
First, the Hardy-Littlewood theorem is true
for $A_{p}$-weights: if $w\in A_{p}$ then
$$
\int_{\bR^{d+1}}\bM^{p}f \,w(t,x)\,dxdt
\leq N\int_{\bR^{d+1}} fw\,dxdt,
$$
where $N$ depends only on $d,p$, and the $A_{p}$-constant $[w]_{p}$
of $w$. In particular, if $r>1$ and $w\in A_{r}$, then
$$
\int_{\bR^{d+1}}\bM^{r}(|f|^{2}) \,w \,dxdt
\leq N\int_{\bR^{d+1}} |f|^{2r}w\,dxdt,
$$

Also, if $w\in A_{r}$, then
there exists $q\in(1,r)$ (close to $r$) and a constant $N$,
depending only on $d$, $r$, and $[w]_{r}$,
such that $w\in A_{q}$. In particular, if $r/\tau\geq q$,
then $w\in A_{r/\tau}$ and
$$
\int_{\bR^{d+1}}\left(\bM(|D^{2}u|^{2\tau})\right) ^{r/\tau}\,w \,dxdt
\leq N\int_{\bR^{d+1}} |D^{2}u|^{2r}w\,dxdt.
$$
This shows how to choose $\tau>1$ depending only
on $d,r$, and $[w]_{r}$.

The final piece of information we need to transform \eqref{eq13.5.05}
is the Fefferman-Stein Theorem which says that for $w\in A_{p}$ 
we have
$$
\int_{\bR^{d+1}}|v|^{p}w\,dxdt\leq N
\int_{\bR^{d+1}}|v^{\#}|^{p}w\,dxdt
$$
if the left-hand side is finite,
where $N$ depends only on $d,p$, and $[w]_{p}$.

By combining all these facts with \eqref{eq13.5.05}
we come to the following.
\begin{corollary}
                                       \label{corollary 11.2.1}
Let $p>2$, $K_{0}\in(1,\infty)$,
$u\in W^{1,2}_{p}(\bR^{d+1})$ and let $w$
be an $A_{p/2}$-weight with $[w]_{p/2}\leq K_{0}$.
 Let $\alpha$ be the constant in
Theorem 4.4 of \cite{DK_10} and $\theta \in (0,1]$.
Then under
Assumption
\ref{assump2} ($\theta$) there exists  constants 
$N$ and $\sigma>1$, depending only on
$d$, $R_{0}$,
$\delta$, $p$, and $K_{0}$, such that,
for any $\kappa\geq4$,
$$
\int_{\bR^{d+1}}|D^{2}_{x''}u|^{p} w\,dxdt\leq N
\kappa^{d+2}\int_{\bR^{d+1}}(|Lu|+|Du|+|u|)^{p} w\,dxdt
$$
\begin{equation}
                                              \label{11.2.4}
+N \kappa^{d+2}\theta^{1/\sigma}
\int_{\bR^{d+1}}|D^{2} u|^{p} w\,dxdt
+N\kappa^{-2\alpha} \int_{\bR^{d+1}}|D^{2}_{x''}u|^{p} w\,dxdt.
\end{equation}

\end{corollary}

To choose $K_{0}$ we use 
Theorem 8.1 of \cite{DK_18_1}
or Theorem 2.5 of \cite{DK_18}  
according to which, for $p,q\in(2,\infty)$
and two functions $f$ and $g$, the inequality
$$
\int_{\bR^{d+1}}|f|^{p/2} w\,dxdt
\leq  \int_{\bR^{d+1}}|g|^{p/2} w\,dxdt,
$$
valid   for any
$A_{p/2}$-weight $w$ with $A_{p/2}$-constant majorated
by $K_{0}(p,q,d)$, for certain constant $K_{0}(p,q,d)$,
implies that
$$
\int_{\bR^{d-2}}\Big(\int_{\bR^{3}}|f|^{p/2}\,dx'dt\Big)^{q/p}
\,dx''\leq N(p,d)
\int_{\bR^{d-2}}\Big(\int_{\bR^{3}}|g|^{p/2}\,dx'dt\Big)^{q/p}
\,dx''.
$$

We fix $p,q\in(2,\infty)$  take $K_{0}=K_{0}(p,q,d)$ 
in Corollary \ref{corollary 11.2.1}
and after that
take $\kappa$ in \eqref{11.2.4} so large that
the last term is absorbed by the left-hand side.

Then  we arrive at

\begin{corollary}
                                       \label{corollary 11.2.2}
Let $p,q>2$ and $u\in W^{1,2}_{p}(\bR^{d+1})$. 
Then under
Assumption
\ref{assump2} ($\theta$) there exists  constants 
$N=N(d,\delta,p,q)$ and $\sigma=\sigma(p,q,d)>1$,  such that 
\begin{equation}
                                              \label{11.2.5}
\|D^{2}_{x''}u\|_{L_{p,q}}\leq N\theta^{1/(p\sigma)}
\|D^{2}u\|_{L_{p,q}}
+ N\||Lu|+|Du|+|u|\|_{L_{p,q}}.
\end{equation}

\end{corollary}

{\bf Proof of Theorem 
\ref{theorem 11.3.1}}.  Let $L_{0}u=\partial_{t}u+a^{ij}D_{ij}u$.
 By Corollary \ref{corollary 2.2}
$$
\|D^{2}u\|_{L_{p,q}}\leq N
\|L_{0}u\|_{L_{p,q}}+N\|D^{2}_{x''}u\|_{L_{p,q}}.
$$
It follows that
$$
\|D^{2}u\|_{L_{p,q}}\leq N(
\|L u\|_{L_{p,q}}+\| Du\|_{L_{p,q}}+\| u\|_{L_{p,q}})
+N\|D^{2}_{x''}u\|_{L_{p,q}}.
$$
This and Corollary \ref{corollary 11.2.2}
shows how to choose small
$\theta=\theta(p,q,d,\delta)>0$ in order to get
\begin{equation}
                                              \label{11.3.4}
\|D^{2}u \|_{L_{p,q}}
\leq N(\| Lu\|_{L_{p,q}}+\| Du\|_{L_{p,q}}+\| u\|_{L_{p,q}}).
\end{equation}
After that to prove assertion (i) of the theorem it suffices to
use Agmon's idea whose implementation
in the mixed norms case the reader can find
in the proof of Lemma 7.2.3 of \cite{Kr_08}.

In light of the method of continuity
to prove assertion (ii) for general $L$ it suffices
to prove it for $L=\partial_{t} +\Delta  $.
Furthermore, assertion (i) implies that it suffices
to show that, if $f\in C^{\infty}_{0}$, then
for $L=\partial_{t} +\Delta  $
there exists a solution of 
\eqref{parabolic} in
$\bR^{d+1} $ of class $W^{1,2}_{p,q}$.

Take such an $f$. Then by the classical theory
there exists a solution $u$ of \eqref{parabolic}
all derivatives of whose are bounded. Next,
take a $\zeta\in C^{\infty}_{0}$ such that $\zeta(0,0)=1$
and for $n=1,2,...$ introduce $\zeta_{n}(t,x)
=\zeta(t/n,x/n)$, $u_{n}=u\zeta_{n}$.
Obviously, $u_{n}\in W^{1,2}_{p,q}$ and, owing to
assertion (i) applied to $u_{n}$,
$$
\lambda\|u \zeta_{n}\|_{L_{p,q}}+
\sqrt{\lambda}\|\zeta_{n}Du \|_{L_{p,q}}+
\|\zeta_{n}D^{2}u,\zeta_{n}\partial_{t}u\|_{L_{p,q}}
\leq N\|f\zeta_{n}\|_{L_{p,q}}
$$
$$
+
N(\sup|u|+\sup|Du|)\big(\|D\zeta_{n}\|_{L_{p,q}}+
\|D^{2}\zeta_{n}\|_{L_{p,q}}
+\|\partial_{t}\zeta_{n}\|_{L_{p,q}}\big),
$$
where the constants $N$ are independent of $n$.
By letting $n\to\infty$ we obtain the desired result.
The theorem is proved.


\end{document}